\documentstyle{article}

\newtheorem{th}{Theorem}
\newtheorem{definition}{Definition}
\newtheorem{lemma}[th]{Lemma}

\title{ Polyhedral representations of discrete differential
manifolds}
\author{Roman R. Zapatrin}
\date{}

\begin{document}
\maketitle
\begin{center}
Department of Mathematics, SPb UEF, Griboyedova 30/32, \\
191023, St-Petersburg, Russia \\
{\em and}
\\ Division of Mathematics, \\
Istituto per la Ricerca di Base, \\
I-86075, Monteroduni, Molise, Italy
\end{center}

\begin{abstract}
Any discrete differential manifold $M$ (finite set endowed
with an algebraic differential calculus) can be represented by
appropriate polyhedron ${\cal P}(M)$. This representation
demonstrates the adequacy of the calculus of discrete differential
manifolds and links this approach with that based on finitary
substitutes of continuous spaces introduced by R.D.Sorkin.
\end{abstract}

\section{ Introduction}

During long time it was the concept of differentiable manifold that
was the arena on which physical theories took place due to its
adequacy to the intuitive feeling what the physical space ought to
be. However, the development of quantum theory gave rise to the
idea to deprive the spacetime of its status of primordial object
and inspired the development of the theories dealing with 'deeper'
entities than spacetime in order to make the latter an observable
in a more general theory. Reasoning about spacetime at very
small scales gave rise to the theories where the manifold was
replaced by a discrete structure. The environment of the present
paper is formed by the following scope of ideas and techniques.

The first belongs to Geroch \cite{geroch} and asserts that
even in classical general relativity the notion of the spacetime
manifold is essentially used only once: to set up the algebra of
smooth functions. We can, instead, start from this (commutative)
algebra as the basic object of the theory. Although. it remains
nothing but a reformulation of the conventional theory and does not
mean that the points (events) are effectively smeared off
\cite{ps}.

If we accept an algebra to be the starting object (called the basic
algebra) we can try to go beyond the class of commutative Banach
algebras representable by functions on manifolds. In particular, we
can assume these algebras to be non-commutative which gives rise to
the non-commutative geometry \cite{connes}.  Another opportunity
looking very attractive from the computational point of view is to
assume the basic algebra to be finite-dimensional and commutative.

At first sight, these objects look very trivial and poor since any
such algebra can be realized by the algebra of functions on a
finite set with the discrete topology. However, even in classical
differential geometry a differentiable manifold is a topological
space plus a differential structure rather than simply a
topological space. Being applied to finite dimensional commutative
algebras this observation gives rise to the notion of discrete
differential calculus and discrete differential manifold
\cite{dmh,dmhv}. This construction will be recalled in sections
\ref{s1} and \ref{s2}.

The important entity in the theory of discrete differential
manifolds is that of the {\em generated topological space} being a
finite or at most countable topological space associated with a
discrete differential manifold (section \ref{s3}). The generated
topological space is intended to play the role of finitary
substitute for the continuous spacetime manifold.

The finitary spacetime substitutes introduced in \cite{sorkin}
serve to simulate the continuous topologgical spaces. These
techniques (reviewed in section \ref{s3}) may be treated as the
generalization of the Regge calculus in general relativity
\cite{regge} for the case when the metric is given up \cite{prg}.
The substitutes are built from continuous manifolds by applying the
coarse-graining procedure yielding finite or at most countable
$T_0$-spaces \cite{icomm}.

So, there are two sources of finitary $T_0$ spaces: the generated
topological spaces of discrete differential manifolds \cite{dmhv}
on one hand, and the results of coarse graining of continuous
topological spaces on the other.

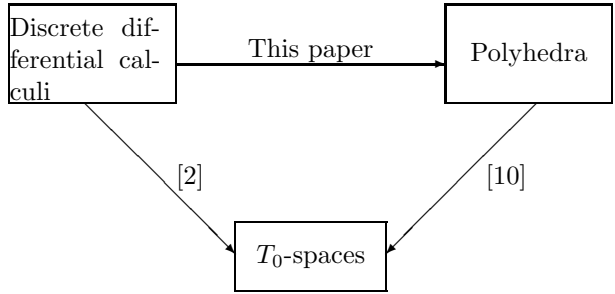
\begin{figure}[hb]
\label{p185}
\unitlength=1mm
\linethickness{0.4pt}
\begin{center}
\begin{picture}(80.00,34.84)
\put(0,25){\framebox(22,13)[cc]{}}
\put(0,30){\parbox{20mm}{Discrete
differential calculi}}
\put(58,25){\framebox(22,13)[cc]{Polyhedra}}
\put(30.00,0.00){\framebox(20.00,9.03)[cc]{$T_0$-spaces}}
\put(22,30.11){\vector(1,0){36}}
\put(10.00,24.95){\vector(1,-1){20.00}}
\put(70.00,24.95){\vector(-1,-1){20.00}}
\put(40,32){\makebox(0,0)[cc]{This paper}}
\put(24,15.05){\makebox(0,0)[cc]{\cite{dmh}}}
\put(66,15.05){\makebox(0,0)[cc]{\cite{icomm}}}
\end{picture}
\end{center}
\caption{The environment of the contents of the paper}
\end{figure}

In this paper the following techniques are suggested. Given a
discrete differential manifold, the polyhedron (being a continuous
topological space) is built and the coarse graining procedure for
it is specified so that the generated topological space of the
discrete differential manifold and the finitary substitute of the
polyhedron are isomorphic $T_0$-spaces (the diagram on Fig.
\ref{p185} is commutative).

Technically, it is done in the following way. A discrete
differential calculus is represented as the quotient of the
universal differential algebra over an appropriate differential
ideal. In section \ref{s2} the structure of differential ideals
is studied. As a result, with any finite dimensional discrete
differential manifold $({\cal M},{\Omega})$ the abstract simplicial complex
${\cal K}={\cal K}({\Omega})$ is associated. In section \ref{s4} the
polyhedral representations of simplicial complexes associated with
discrete differential manifolds are considered, and the
coarse-graining procedure \cite{sorkin,icomm} described in section
\ref{s3} is specified in order to produce the finitary substitutes
isomorphic to the topological spaces generated by the discrete
differential manifolds.

\section{ Basic notions and results} \label{s1}

Let ${\cal M}$ be a finite set. An algebraic {\sc differential
calculus} \cite{dmh} on ${\cal M}$ is an extension of the algebra ${\cal A} =
{\sf Fun}({\cal M},{\bf C})$ (called {\sc basic algebra} of all ${\bf C}$-valued
functions on ${\cal M}$ to a graded differential algebra $({\Omega}, d)$:

\begin{equation}\label{f187o}
{\Omega} = {\Omega}^0 \oplus {\Omega}^1 \oplus \ldots \oplus {\Omega}^r \oplus \ldots
\end{equation}
with ${\Omega}^0={\cal A}$, $d:{\Omega}^r\to{\Omega}^{r+1}$ and $d{\Omega}^r$ being the
generating set for the next ${\cal A}$-bimodule ${\Omega}^{r+1}$. If
${\Omega}^r\neq 0$ for some $r$ and ${\Omega}^s=0$ for every $s>r$, then
the calculus ${\Omega}$ is said to have the {\sc dimension} $r$. The
universal object ${\tilde{\Omega}}={\tilde{\Omega}}({\cal M})$ of such type is the {\sc
differential envelope} \cite{kastler} of ${\cal A}$. The set ${\cal M}$ is
assumed to be finite, and ${\tilde{\Omega}}$ can be explicitly described by
setting its {\sc natural basis}:
\begin{equation}\label{f187b}
\begin{array}{rcl}
{\tilde{\Omega}}^0 &=& {\sf span}\{e_i\mid i\in{\cal M}\} \cr
{\tilde{\Omega}}^1 &=& {\sf span}\{e_{ik}\mid i,k\in{\cal M}\hbox{ , }i\neq k\}\cr
\ldots &\ldots& \ldots \cr
{\tilde{\Omega}}^r &=& {\sf span}\{e_{i_0,i_1,\ldots,i_r}\mid
i_0,i_1,\ldots,i_r\in{\cal M}\hbox{ ; }\forall s=1,\ldots,r\quad
i_{s-1}\neq i_s\}
\end{array}
\end{equation}
where $e_i:{\cal M}\to{\bf C}$ is defined as $e_i(k)=\delta_{ik}$. Then each
${\tilde{\Omega}}^r$ is ${\cal A}$-bimodule:
\begin{equation}\label{f187a}
e_p\cdot e_{i_0,i_1,\ldots,i_r}\cdot e_q =
\delta_{pi_0}\delta_{i_rq}e_{i_0,i_1,\ldots,i_r}
\end{equation}
and the graded product ${\tilde{\Omega}}^r\times{\tilde{\Omega}}^s\to{\tilde{\Omega}}^{r+s}$ is defined
as:
\begin{equation}\label{f187p}
e_{i_0,i_1,\ldots,i_r}e_{j_0,j_1,\ldots,j_s} =
\delta_{i_rj_0}e_{i_0,i_1,\ldots,i_r,j_1,\ldots,j_s}
\end{equation}
with the operator $d:{\tilde{\Omega}}^r\to{\tilde{\Omega}}^{r+1}$ having the form:
\begin{equation}\label{f187d}
de_{i_0,i_1,\ldots,i_r} =
\sum_{s=0}^{r}(-1)^s
\sum_{k\neq i_s,i_{s+1}}e_{i_0,\ldots,i_{s-1},k,i_s,\ldots,i_r}
\end{equation}

The graded algebra ${\tilde{\Omega}}$ is universal in the sense that any
particular differential calculus ${\Omega}$ (\ref{f187o}) can be
covered by an epimorphism $\pi:{\tilde{\Omega}}\to{\Omega}$ of graded
differential algebras over the basic algebra ${\cal A}$. The kernel
of this mapping $\pi$ is said to be {\sc differential ideal} in
${\tilde{\Omega}}$. So, every differential calculus ${\Omega}$ can be unambiguously
characterized by appropriate differential ideal ${\cal I}({\Omega})$ in
${\tilde{\Omega}}$. Let us dwell on this issue in more detail.

\begin{definition}
A linear subspace ${\cal I}\subseteq{\tilde{\Omega}}$ is called {\sc
differential ideal} in ${\tilde{\Omega}}$ if:
\begin{equation}\label{f188}
\begin{array}{rcl}
{\tilde{\Omega}}{\cal I}{\tilde{\Omega}} &\subseteq& {\cal I} \cr
d{\cal I} &\subseteq& {\cal I}
\end{array}
\end{equation}
\end{definition}

Denote by ${\Omega}({\cal I})={\tilde{\Omega}}/{\cal I}$ the differential calculus
(\ref{f187o}) induced by the differential ideal ${\cal I}$. The
decomposition of the graded algebra ${\tilde{\Omega}}$ gives rise to the
decomposition of ${\cal I}$:
\begin{equation}\label{f189i}
{\cal I} = {\cal I}^0 \oplus \ldots \oplus {\cal I}^r \oplus \ldots
\end{equation}
Then the definition of the differential calculus ${\Omega}({\cal I})$ is
reformulated in terms of ${\cal I}$ as follows:
\[ \dim{\Omega}({\cal I})  = r\quad
\stackrel{\rm def}{\Leftrightarrow}\quad
{\cal I}^r\neq{\tilde{\Omega}}^r \hbox{ and }
\forall s>r\quad {\cal I}^s={\tilde{\Omega}}^s
\]

\paragraph{Assumption.} We shall consider only the differential
ideals with ${\cal I}^0=0$.

In the sequel, to introduce combinatorial structures confine
ourselves by {\sc basic differential ideals}, namely, the ideals
spanned on the elements of the natural basis (\ref{f187b}).

\begin{lemma}\label{l189} Let ${\cal I}$ be a basic differential ideal,
$\alpha = i_0,\ldots,i_r$ and $e_\alpha\in{\cal I}$. Let $\beta$ be such
sequence of elements of ${\cal M}$ that $\alpha$ is its subsequence
$\alpha\subseteq\beta$. Then $e_\beta\in{\cal I}$.
\end{lemma}

\noindent {\bf Proof.} Use the induction over $\delta = \mid\beta\mid -
\mid\alpha\mid$ (the difference of the lengths of the sequences).
Suppose $\delta=1$, that means $\beta = i_0 ,\ldots, i_{t-1},q,i_t,
\ldots, i_r$. Due to (\ref{f188}), $de_\alpha\in {\cal I}$ and equals to
(\ref{f187d}), being the sum of the basis elements of ${\cal I}$
(assumed to be the basic ideal). Therefore each summand of
(\ref{f187d}), in particular, $e_\beta$, is the element of ${\cal I}$.
Now suppose the result is valid for all $\beta$ such that
$\mid\beta\mid - \mid\alpha\mid = s$, and let $\beta$ be such
supersequence of $\alpha$ that $\mid\beta\mid = \mid\alpha\mid + s
+ 1$. Consider a sequence $\beta'$ obtained from $\beta$ by
deleting an element which does not belong to the sequence $\alpha$.
Then $e_{\beta'}\in{\cal I}$ (being the supersequence of $\alpha$ of
the length $\mid\alpha\mid + s$), and therefore $e_\beta\in{\cal I}$
since it occurs in the decomposition (\ref{f187d}) for
$de_{\beta'}$.  \medskip

\paragraph{Corollary.} Let ${\cal I}$ is a basic differential ideal, and
$\alpha,\beta$ are such sequences of elements of ${\cal M}$ that
$\alpha\subseteq \beta$. Then $e_\beta\not\in{\cal I} \Rightarrow
e_\alpha\not\in{\cal I}$.

It is suitable to impose the scalar product on ${\tilde{\Omega}}$ by assuming
the natural basis (\ref{f187b}) to be orthonormal:
\begin{equation}\label{f61s}
(e_\alpha,e_\beta) = \left\lbrace \begin{array}{rcl}
1 &,& \hbox{if}\quad \alpha=\beta \cr
0 &&  \hbox{otherwise}
\end{array}\right.
\end{equation}

\begin{definition} \label{df60} A {\sc basic differential
calculus} is the quotient ${\Omega} = {\Omega}({\cal I}) = {\tilde{\Omega}}/{\cal I}$ of the
differential envelope over a basic differential ideal ${\cal I}$.
\end{definition}

Using the scalar product (\ref{f61s}) the following description of
basic differential calculi was suggested in \cite{dmh}. The
quotient space ${\Omega} = {\tilde{\Omega}}/{\cal I}$ can be identified with the subspace
of ${\tilde{\Omega}}$ being the orthocomplement to ${\cal I}$, denote it with the
same symbol ${\Omega}$. Thus

\[ {\Omega} = {\sf span}\{e_{i_0,\ldots,e_r}\mid \quad
e_{i_0,\ldots,e_r}\not\in{\cal I} \} \]

\noindent or, in other words, the basic differential forms
constituting the ideal ${\cal I}$ are said to be {\em vanishing} in the
differential calculus ${\Omega}$:  \begin{equation}\label{f61n} e_{i_0,\ldots,i_r}\in {\cal I}
\quad\Leftrightarrow\quad e_{i_0,\ldots,i_r}=0
e_{i_0,\ldots,i_r}\not\in {\cal I} \quad\Leftrightarrow\quad
e_{i_0,\ldots,i_r}\neq0
\end{equation}

So, the Lemma \ref{l189} can be reformulated as follows:
\begin{equation}\label{f189l}
\forall \alpha \subseteq \beta \qquad e_\beta\neq 0 {\quad\Rightarrow\quad}
e_\alpha \neq 0
\end{equation}

\section{Discrete differential manifolds}\label{s2}

In classical differential geometry a differentiable manifold is a
topological space equipped with a differential structure. Its
finitary counterpart looks as follows.

\begin{definition} A {\sc discrete differential manifold} is a
couple $({\cal M},{\Omega})$ where ${\cal M}$ is a finite set and ${\Omega}$ is a basic
differential calculus (see Definition \ref{df60}) over the
functional algebra ${\cal A}={\sf Fun}({\cal M})$. The {\sc dimension} of the
discrete differential manifold is the greatest grade of its
nonvanishing differential forms:  \[ \dim {\cal M} = \max\{r\mid
{\Omega}^r\neq 0\} \] If this number does not exist, the discrete
differential manifold is said to be infinite dimensional:
\begin{equation}\label{f63d} \dim {\cal M} = \infty \stackrel{\rm
def}{\Leftrightarrow}\quad \forall r \quad {\Omega}^r\neq 0 \end{equation}
\end{definition}

Taking into  account the notation (\ref{f61n}), discrete
differential manifolds are unambiguously determined by the
collection of their nonvanishing basic differential forms:
\[
M = ({\cal M},{\Omega}) = ({\cal M},{\cal K})
\]
where ${\cal K} = {\cal K}(M)$ is the collection of the sequences $\alpha$ of
elements of the set ${\cal M}$ such that
\begin{equation}\label{f63k}
\alpha\in {\cal K}
\stackrel{\rm def}{\Leftrightarrow}\quad e_\alpha \neq 0
\end{equation}

The set of sequences ${\cal K}$ can be decomposed in the following way:
\[
\begin{array}{rcl}
{\cal K} &=& {\cal K}^1 \cup {\cal K}^2 \cup \ldots \cup {\cal K}^r \cup \ldots \cr
\alpha\in {\cal K}^r &{\quad\Leftrightarrow\quad}& \alpha\in {\cal K} \quad\hbox{and}\quad \vert
\alpha \vert = r
\end{array}
\]

\begin{lemma} \label{l62} If $(i_0,\ldots,i_r)\in {\cal K}^r$ then for
any $s$ such that $0\le s\le r$ we have
\( (i_0,\ldots,\hat{i_s},\ldots,i_r)\in {\cal K}^{r-1} \)
\end{lemma}

\noindent {\bf Proof.} Follows from the corollary from Lemma \ref{l189}.  \medskip

The component ${\cal K}^1$ of the collection ${\cal K}$ gives rise to a
binary relation on the set ${\cal M}$, denote it $\preceq$:
\begin{equation}\label{f63p}
i\preceq j
\stackrel{\rm def}{\Leftrightarrow}\quad (i,j)\in {\cal K}^1 {\quad\Leftrightarrow\quad}
e_ij\not\in {\cal I}
\end{equation}
Since no $e_{ii}$ is the basic form, $e_{ii}\not\in {\cal I}$, hence the
relation $\preceq$ on ${\cal M}$ is {\em reflexive}: \( \forall i\in {\cal M}
\quad i\preceq i \). It follows immediately from Lemma \ref{l62}
that

\begin{equation}\label{f63}
(i_0,\ldots,i_r)\in {\cal K}^r {\quad\Rightarrow\quad} \forall s,t \quad 0\le s \le t
\le r \quad i_s\preceq i_t
\end{equation}

\begin{definition} \label{df63} A discrete differential manifold
$M$ is called {\sc network manifold} if its differential structure
is completely defined on the level of its 1-forms, that is the
reverse of the implication (\ref{f63}) holds:
\begin{equation}\label{f63n}
{\cal K}^r = \{ (i_0,\ldots,i_r)\mid
\forall s,t \quad 0\le s \le t \le r \quad i_s\preceq i_t \}
\end{equation}
\end{definition}

\begin{lemma}\label{l64} Let $M=({\cal M},{\Omega})$ be a discrete
differential manifold such that the set ${\cal M}$ is finite. If $M$ is
infinite dimensional, then the relation $\preceq$ on the set ${\cal M}$
is not antisymmetric:
\[
\dim M = \infty {\quad\Rightarrow\quad}
\exists i,j\in {\cal M}, i\neq j:\quad i\preceq j, j\preceq i
\]

If $M$ is the network manifold, the above implication holds in both
directions.
\end{lemma}

\noindent {\bf Proof.} Let $m={\sf card} {\cal M}$. Consider a number $r>m$ and let
$e_{i_0,\ldots,i_r}\neq 0$: it exists due to (\ref{f63d}). Then at
least one element $i\in {\cal M}$ occurs at least twice in the string
$(i_0,\ldots,i_r)$, that is for some $s,t$ such that $s>t+1$ we
have $i=i_s=i_t$. Denote $j=i_{s+1}$, then $j\neq i$ by virtue of
(\ref{f187b}), and according to (\ref{f63}) $i\preceq j$ and
$j\preceq i$.

Now let $M$ be a network manifold. Let $i\neq j$ and $i\preceq
j\preceq i$. For any $r>0$ consider the sequence $(ij\ldots ij)$ of
length $2r$. Then, according to (\ref{f63n}) $e_{ij\ldots ij}\neq
0$, thus $\dim M = \infty$.  \medskip

\paragraph{Corollary.} Let $\dim M < \infty$. If $e_\alpha\neq 0$
in $M$, then all the elements of the string $\alpha$ are different.
Therefore
\begin{equation}\label{f65c}
{\sf card}\{i_0,\ldots,i_r\} = {\sf length}(i_0,\ldots,i_r)
\end{equation}

>From now on confine ourselves by finite dimensional discrete
differential manifolds and recall further definitions.

\begin{definition} A set $K$ with a relation $\preceq$ on it is
called {\sc filly ordered} if the relation $\preceq$ has the
following properties: for any $i,j,k\in K$
\begin{equation}\label{f65f}
\begin{array}{rcll}
 &i\preceq i& &
\hbox{\rm (reflexivity)} \cr
i\preceq j, j\preceq i &\Rightarrow& i=j &
\hbox{\rm (antisymmetry)} \cr i\preceq j, j\preceq k &\Rightarrow&
i\preceq k & \hbox{\rm (transitivity)} \cr i\not\preceq j
&\Rightarrow& j\preceq i & \hbox{\rm (linearity)} \end{array} \end{equation}
\end{definition}

Note that any subset of a fully ordered set is fully ordered as
well. The following lemma shows the relevance of full orders in the
theory of discrete differential manifolds.

\begin{lemma} \label{l65} Let $M=({\cal M},{\Omega})$ be a finite dimensional
discrete differential manifold. Consider a {\em sequence} $\alpha$.
Then the following implication holds:
\begin{equation}\label{f65}
e_\alpha\neq 0 {\quad\Rightarrow\quad} (\alpha,\preceq)\quad \hbox{\rm is the full
order}
\end{equation}
Moreover, if $M$ is a network manifold, the implication
(\ref{f65}) holds in both directions.
\end{lemma}

\noindent {\bf Proof.} The implication (\ref{f65}) follows from (\ref{f63}) and the
definition (\ref{f65f}). Now let $M$ be a network manifold, and
$\alpha=\{i_0,\ldots,i_r\}$ be a fully ordered (with respect to
$\preceq$) subset of ${\cal M}$. Arrange the elements of the {\em set}
$\alpha$ to form the {\em sequence} $\alpha$ such that $i_0\preceq
i_1\preceq \ldots \preceq i_r$. Then the transitivity of the
relation $\preceq$ implies that for all $s,t$ such that $0\le s<t
\le r$ we have $e_{i_si_t}\neq 0$, therefore $e_\alpha\neq 0$ by
virtue of (\ref{f63n}).  \medskip

So, from the combinatorial point of view the discrete differential
manifolds are characterized in the following way. We have a set
${\cal M}$ with a reflexive antisymmetric (but not generally transitive)
relation $\preceq$ on it. Then we select a family ${\cal K}$ of subsets
of ${\cal M}$ such that
\begin{itemize}
\item any element $\alpha\in {\cal K}$ is fully ordered by the relation
$\preceq$
\item ${\cal K}$ is hereditary: $\alpha\in {\cal K}, \beta\subseteq \alpha
\Rightarrow \beta\in {\cal K}$
\item ${\cal K}$ contains all singletons (since $\preceq$ is reflexive):
$\forall i\in {\cal M} \quad \{i\}\in {\cal K}$
\end{itemize}

Then we can build the discrete differential calculus ${\Omega}$ on ${\cal A}
= {\sf Fun}({\cal M})$ by putting
\[ {\Omega} = {\Omega}({\cal K}) = {\sf span}\{e_\alpha\mid \alpha \in {\cal K} \} \]
In particular, when $M$ is a finite dimensional network discrete
differential manifold, the appropriate family ${\cal K}({\cal M})$ is the
collection of {\em all} $\preceq$-fully ordered subsets of ${\cal M}$.

To conclude the section, note that the properties of the family
${\cal K}$ listed above coincide with the definition of ordered
simplicial complex. We shall return to this issue in section
\ref{s4}.

\section{ Two sources of finite topological spaces} \label{s3}

The first source described in this section is the procedure
manufacturing finite $T_0$-topological spaces (called {\em
generated topological spaces} \cite{dmh}) from discrete
differential manifolds. The second source is the coarse graining
procedure applicable to arbitrary topological spaces which also
yields the finite $T_0$-spaces \cite{sorkin,icomm}. Begin with the
first source.

\paragraph{Generated topological spaces.} Let $M=({\cal M},{\Omega})$ be a
discrete differential manifold and ${\cal K} = {\cal K}({\cal M})$ be the
collection of its nonvanishing basic forms\footnote{ In \cite{dmhv}
the set ${\cal K}$ is denoted by $\hat{{\cal M}}$} (\ref{f63k}).  Define the
topology $\tau$ on the set ${\cal K}$ by setting its prebase of open
sets $\{U_\alpha\}, \alpha\in {\cal K}$:  \[ U_\alpha=\{\beta\in {\cal K}\mid
\alpha\subseteq \beta\} \]

\begin{definition}\label{df67} The topological space $({\cal K},\tau)$
defined above is called the {\sc generated topological space}
$({\cal K}, \tau)$ of the discrete differential manifold $M$.
\end{definition}

It was shown in \cite{dmhv} that the topology of any generated
topological space is always $T_0$ (that is, for each pair of
distinct points of ${\cal K}$ there is an open set containing one point
but not the other).

\paragraph{Example.} Let ${\cal M} = \{1,2,3\}$. Define the relation
$\preceq$ on ${\cal M}$ as $1\prec 2\prec 3\prec 1$ (not being
transitive), the graph of the relation is on Fig. \ref{p69}a. Let
$M$ be the {\em network} manifold induced by the relation $\preceq$
on ${\cal M}$, then
\[ {\Omega} = {\sf span}\{e_{12},e_{23},e_{31}\} \]
and hence $\dim M = 1$. The generated space ${\cal K}$ is:
\[ {\cal K} = \{1,2,3,12,23,31\} \]
and the topology $\tau$ on ${\cal K}$ is depicted on Fig. \ref{p69}b in
terms of the appropriate Hasse graph (see \cite{dmh,sorkin,icomm}
for details).

\begin{figure}[hb]
\label{p69}
\unitlength=1mm
\linethickness{0.4pt}
\begin{center}
\begin{picture}(100,40)
\put(10,15){\circle*{2}}
\put(10,12){\makebox(0,0)[cc]{1}}
\put(10,15){\vector(1,0){28}}
\put(40,15){\circle*{2}}
\put(40,12){\makebox(0,0)[cc]{2}}
\put(40,15){\vector(-1,1){13}}
\put(25,30){\circle*{2}}
\put(25,34){\makebox(0,0)[cc]{3}}
\put(25,30){\vector(-1,-1){13}}
\put(24,0){\makebox(0,0)[cc]{a).}}

\put(70,15){\line(0,1){15}}
\put(70,15){\line(1,1){15}}
\put(85,15){\line(0,1){15}}
\put(85,15){\line(1,1){15}}
\put(100,15){\line(0,1){15}}
\put(100,15){\line(-2,1){30}}

\put(70,15){\circle*{2}}
\put(71,12){\makebox(0,0)[cc]{12}}
\put(85,15){\circle*{2}}
\put(86,12){\makebox(0,0)[cc]{23}}
\put(100,15){\circle*{2}}
\put(101,12){\makebox(0,0)[cc]{31}}
\put(70,30){\circle*{2}}
\put(70,34){\makebox(0,0)[cc]{1}}
\put(85,30){\circle*{2}}
\put(85,34){\makebox(0,0)[cc]{2}}
\put(100,30){\circle*{2}}
\put(100,34){\makebox(0,0)[cc]{3}}
\put(84,0){\makebox(0,0)[cc]{b).}}
\end{picture}
\end{center}
\caption{a). The graph of the relation $\preceq$. b). The Hasse
graph of the topological space ${\cal K}$. A set $A\subseteq {\cal K}$ is
open iff with every its element $a\in A$ it contains all elements
lying below $a$ and linked with it.}
\end{figure}

\paragraph{Finitary substitutes.} Let $V$ be a topological space
and let ${\cal T} = \{V_\alpha\}$ be its finite open covering: $V =
\cup V_\alpha$. Define the new topology $\tau$ on $V$ as that
generated by the collection ${\cal T}$ considered prebase of open sets.
The set $(V,\tau)$ is in general not even $T_0$ space therefore the
theorem of the uniqueness of limits of sequences may not hold.
Define the relation denoted $\to$ on the set $V$ as follows:
\[ x\to y {\quad\Leftrightarrow\quad} y = \lim_{\tau}\{x,x,\ldots,x,\ldots\} \]
where $\lim_{\tau}$ denotes the limit with respect to the
topology $\tau$ or, in a more transparent form
\[ x\to y {\quad\Leftrightarrow\quad}
(\forall \alpha \quad y\in V_\alpha \Rightarrow x\in V_\alpha ) \]

In general, the relation $\to$ is a preorder on the set $V$ hence
we can consider its quotient ${\cal K} = V/\sim$ with respect to the
equivalence $\sim$:
\begin{equation}\label{f70} x\sim y {\quad\Leftrightarrow\quad} x\to y, y\to x \end{equation}
As a result, the set ${\cal K}$ is partially ordered by the relation
$\to$ and the topology $\tau$ induced on ${\cal K}$ as the quotient set
is already the $T_0$-topology. The detailed account of this
procedure can be found in \cite{sorkin,prg}.

\paragraph{Example.} Let $V$ be a circle, $V=e^{i\phi}$. Consider
the covering ${\cal T} = \{V_\alpha,V_\beta,V_\gamma\}$:
\[
\begin{array}{rcl}
V_\alpha &=& \{e^{i\phi}\mid -\pi/2 <\phi <+\pi \} \cr
V_\beta &=& \{e^{i\phi}\mid \pi/2 <\phi <3\pi/4 \} \cr
V_\gamma &=& \{e^{i\phi}\mid -\pi <\phi <+\pi/4 \}
\end{array}
\]
Then the equivalence classes (\ref{f70}) (that is, the elements of
${\cal K}$ are:
\[
\begin{array}{rcl}
1 &=& \{e^{i\phi}\mid -\pi/2 <\phi <+\pi/4 \} \cr
12 &=& \{e^{i\phi}\mid \pi/4 <\phi <\pi/2 \} \cr
2 &=& \{e^{i\phi}\mid \pi/2 <\phi <\pi \} \cr
23 &=& \{e^{i\pi}\} \cr
3 &=& \{e^{i\phi}\mid -\pi <\phi <-\pi/2 \} \cr
31 &=& \{e^{-i\pi/2}\}
\end{array}
\]
and the induced partial order has the Hasse diagram the same as
depicted on Fig. \ref{p69}b.

It will be shown in section \ref{s4} how, starting from an
arbitrary finite dimensional discrete differential manifold to
build a continuous metrical space (namely, a polyhedron) and
specify its open covering so that the resulting $T_0$ spaces would
be the same. To do it, we have to introduce one more technical
issue.

\paragraph{Simplicial coarse graining of polyhedra.} Let ${\cal P}$ be a
simplicial complex and ${\vert{\cal P}\vert}$ be its realization by a polyhedron in
a Euclidean space ${\cal E}$. That means that ${\vert{\cal P}\vert}$ is the union of
well-positioned geometrical simplices $\alpha$ in ${\cal E}$.
Initially, ${\vert{\cal P}\vert}$ is the metrical space being the subset of the
space ${\cal E}$ with the standard Euclidean metric and having the
topology associated with this metric. For every point $x\in {\vert{\cal P}\vert}$
of the polyhedron consider its {\em star}:
\[ {\sf St}(x) = \cup\{\alpha \mid x\in \alpha\} \]
and then define the neighborhood of the point $x$ as the interior
of ${\sf St}(x)$ in ${\vert{\cal P}\vert}$:
\[ V_x = {\sf Int}_{\vert{\cal P}\vert} {\sf St}(x) \]
(note that for some $x\neq y$ the neighborhoods may coincide).
Evidently ${\cal T} = \{V_x\mid x\in {\vert{\cal P}\vert}\}$ is the open covering of
${\vert{\cal P}\vert}$. Moreover, infinite as the set ${\vert{\cal P}\vert}$ is, the covering
${\cal T}$ is finite. The elements of ${\cal T}$ are in 1-1 correspondence
with the simplices of ${\cal P}$.  \[ (\forall \alpha \in {\cal P} \quad x\in
\alpha \Leftrightarrow y\in \alpha)\quad \hbox{iff} \quad V_x=v_y
\]

Now consider the topology $\tau$ on ${\vert{\cal P}\vert}$ induced by the covering
${\cal T}$ thought of as open prebase, and the appropriate
$T_o$-quotient. For every $V_\alpha, V_\beta \in {\cal T}$ the
intersection $V_\alpha\cap V_\beta =V_{\alpha\cap \beta}$ is either
$\emptyset$ or the element of ${\cal T}$, therefore ${\cal T}$ is the base
(rather than prebase) of the topology $\tau$.

\begin{lemma} There is the 1-1 correspondence between the simplices
of ${\cal P}$ and the points of the quotient space:
\[ {\cal P} = {\vert{\cal P}\vert}/\sim \]
\end{lemma}

\noindent {\bf Proof.} Associate with every simplex $\alpha \in {\cal P}$ its 'local
interior' $I(\alpha)$:
\[ I(\alpha) = \{(\mu^{\alpha}_0,\ldots,\mu^{\alpha}_r)\mid \forall
i = 1,\ldots ,r \quad 0<\mu^{\alpha}_i<1\} \]
where $r$ is the dimension of the simplex $\alpha$ and
$\mu^{\alpha}_i$ are its baricentric coordinates. For the vertices
$v\in {\cal P}$ put:
\[ I(v) = \{v\} \]
The collection $\{I(\alpha)\mid \alpha \in {\cal P}\}$ is the partition
of the polyhedron ${\vert{\cal P}\vert}$. For any $\alpha$, $x,y \in I(\alpha)$
implies $V_x = V_y$, and vice versa $V_x = V_y$ implies $\exists
\alpha x,y\in I(\alpha)$.

So, we may conclude that the resulting quotient $T_0$-space is the
simplicial complex ${\cal P}$ itself (that is, the points of the
finitary substitute are the simplices of ${\vert{\cal P}\vert}$). The open base of
the topology $\tau$ on ${\cal P}$ is the collection of the stars of all
simplices of ${\cal P}$.  \medskip

The appropriate open covering of the polyhedron ${\vert{\cal P}\vert}$ is called
{\em simplicial} since it does not depend on particular realization
${\vert{\cal P}\vert}$ of the complex ${\cal P}$.

\paragraph{Example.} Let the polyhedron ${\vert{\cal P}\vert}$ be the triangle
without interior whose vertices are labelled by 1,2,3. then the
appropriate simplicial complex is ${\cal P} = \{1,2,3,12,13,23\}$ (for
brevity I denote $1 = \{1\}, 12 = \{1,2\}$ and so on).  The
simplicial covering $\tau$ consists of 6 sets:
\[
\begin{array}{rclcrcl}
V_{12} &=& (1,2) &\hbox{        }& V_1 &=& (1,2)\cup(1,3) \cr
V_{23} &=& (2,3) &\hbox{        }& V_2 &=& (1,2)\cup(2,3) \cr
V_{13} &=& (1,3) &\hbox{        }& V_3 &=& (1,3)\cup(2,3)
\end{array}
\]
where $(\cdot,\cdot)$ denotes the open interval between appropriate
vertices. Then the finitary substitute induced by the simplicial
covering ${\cal T}$ is finite topological space depicted on the Fig.
\ref{p69}b.

\section{ Polyhedral representations of discrete differential
manifolds and the correspondence theorem} \label{s4}

In this section the main result of the paper is formulated, namely
the transition from discrete differential manifolds to polyhedra
associated with the same $T_0$ spaces is described. Let
$M=({\cal M},{\Omega})$ be a finite dimensional discrete differential
manifold.

\begin{lemma}\label{l78} Let $\{i_0,\ldots, i_r\}$ be a subset of
${\cal M}$.  Then there exists at most one basic differential form
$e_\alpha$ such that \begin{itemize} \item $\alpha$ is a
permutation of the elements $\{i_0,\ldots, i_r\}$.  \item $e_\alpha
\neq 0$ \end{itemize} \end{lemma}

\noindent {\bf Proof.} Suppose there are $e_\alpha, e_\beta$ such that $\beta$ is
obtained from $\alpha$ by a nontrivial permutation $\sigma$. Let
$\alpha = (i_0,\ldots, i_r)$. Since $\sigma\neq {\sf id}$ there
exists a pair of distinct elements $i,j\in {\cal M}$ such that $i$
precedes $j$ in $\alpha$ and $j$ precedes $i$ in $\beta$. Thus it
follows from (\ref{f63}) that $i\preceq j$ and $j\preceq i$ which
contradicts with the assumption $\dim M<\infty$ (lemma \ref{l64}).
 \medskip

With any finite dimensional discrete differential manifold
$M=({\cal M},{\Omega})$ its generated topological space ${\cal K}(M)$ (definition
\ref{df67}) may be considered as the collection of {\em subsets}
(rather than sequences) of the set ${\cal M}$: due to lemma \ref{l78} we
can forget about the order and different ordered sets will become
different sets (note that this does not work when $\dim M=\infty$
!)

Recalling the properties of ${\cal K}(M)$ established in the end of
section \ref{s2} we see that ${\cal K}({\cal M})$ is exactly the simplicial
complex with the set of vertices ${\cal M}$.

\begin{definition} A {\sc polyhedral representation} of a discrete
differential manifold M is the polyhedron $\vert {\cal K}\vert$ being a
geometrical realization of the simplicial complex ${\cal K}$.
\end{definition}

Now the main result of the paper can be formulated as the following
{\em correspondence theorem}:

\begin{th} Let $M=({\cal M},{\Omega})$ be a finite dimensional discrete
differential manifold, ${\cal K}(M)$ be its polyhedral realization. The
following two $T_0$-spaces are isomorphic:
\begin{itemize}
\item
The finitary substitute $\vert {\cal K}\vert/{\cal T}$ with respect to the
simplicial covering ${\cal T}$ of $\vert {\cal K}\vert$.  \item The
generated topological space ${\cal K}$ of the discrete differential
manifold $M$.  \end{itemize} \end{th}

\noindent {\bf Proof.} In fact, everything is already proved. The points of both
$T_0$-spaces are in 1-1 correspondence with the elements of the
complex ${\cal K}(M)$. The partial orders on both finite sets ${\cal K}$
and $\vert {\cal K}\vert/{\cal T}$ are the same being the set inclusion of
simplices. Thus ${\cal K}$ and $\vert {\cal K}\vert/{\cal T}$ are homeomorphic.
 \medskip

\paragraph{Remark.} In \cite{dmhv} the question which finite
$T_0$-spaces are generated by a discrete differential manifold
remained open. In this paper these $T_0$ spaces are characterized.
Moreover, it is seen that simplicial complexes and related
structures such as polyhedra or simplicial spaces are more adapted
to be the topological realizations for discrete differential
manifolds rather than finitary topological spaces of general form.

\medskip
\paragraph{Acknowledgments.} The author is grateful to A.Dimakis
and the participants of the Efroimsky seminar (St.Petersburg),
particularly, G.N. Parfionov, for helpful discussions and remarks.


\begin{thebibliography}{99}
\bibitem{connes} Connes A., Noncommutative differential geometry,
Hermann, Paris, 1989

\bibitem{dmh} Dimakis A., F. M\"uller-Hoissen,
{\it Discrete differential calculus, graphs, topologies and gauge
theory},
Journal of Mathematical Physics,
{\bf 35},
6703,
(1994)

\bibitem{dmhv} Dimakis A., F. M\"uller-Hoissen, F. Vanderseypen,
{\it Discrete differential manifolds and dynamics on networks},
Journal of Mathematical Physics,
{\bf 36},
3771,
(1995)
(eprint hep-th/9408114)

\bibitem{geroch} Geroch, R.
{\it Einstein Algebras},
Communications in Mathematical Physics,
{\bf 26},
271,
(1972)

\bibitem{kastler} Kastler, D.,
{\it Cyclic cohomology within the differential envelope},
Hermann, Paris,
1988

\bibitem{ps} G.N.Parfionov, R.R.Zapatrin,
{\it Pointless Spaces in General Relativity},
International Journal of Theoretical Physics,
{\bf 34}, 737,
1995 (eprint gr-qc/9503048)

\bibitem{regge} T. Regge,
{\it General relativity without coordinates},
Nuovo Cimento,
{\bf 19},
568,
1961

\bibitem{sorkin} Sorkin, R.D.,
{\it Finitary substitutes for continuous topology},
International Journal of Theoretical Physics,
{\bf 30},
923,
(1991)

\bibitem{prg} Zapatrin, R.R.,
{\it Pre-Regge Calculus: Topology Via Logic},
International Journal of Theoretical Physics,
{\bf 32},
779,
(1993)

\bibitem{icomm} Zapatrin, R.R.,
{\it Matrix models for spacetime topodynamics},
In: Proceedings of the ICOMM'95 (Vienna, June 3-6, 1995),
W. Kainz, Ed., TUV (1995),
1-19
(eprint gr-qc/9503066)

\end{thebibliography}
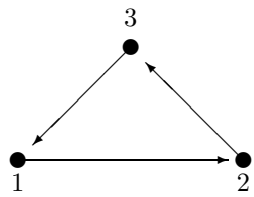
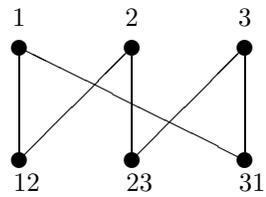
\end{document}